\documentclass[11pt]{article}
\usepackage{amsfonts}
\usepackage{amsmath}
\usepackage{amssymb}
\usepackage{graphicx}
\newfont{\blb}{msbm10 scaled\magstep1}
\newfont{\comp}{cmr12 scaled\magstep1}
\newfont{\compb}{cmr10 scaled\magstep2}
\newtheorem{theorem}{Theorem}[section]
\newtheorem{lemma}[theorem]{Lemma}

\newtheorem{corollary}[theorem]{Corollary}

\newtheorem{claim}[theorem]{Claim}

\newenvironment{proof}{{\bf Proof.}}{\hfill{ }\vrule height10pt width5pt depth1pt}

\date{}

\begin{document}

\author{Benny Sudakov \thanks{
Department of Mathematics, Princeton University, Princeton, NJ 08544, and
Institute for Advanced Study, Princeton. E-mail:
bsudakov@math.princeton.edu.
Research supported in part by NSF CAREER award DMS-0546523, NSF grant
DMS-0355497, USA-Israeli BSF grant, Alfred P. Sloan fellowship, and
the State of New Jersey.}
\and Jacques Verstra\"{e}te
\thanks{Department of Combinatorics and Optimization, University of
Waterloo, Waterloo, ON N2L 3G1, Canada. E-mail:
jverstra@math.uwaterloo.ca}}

\title{Cycle lengths in sparse graphs}

\maketitle

\vspace{-0.3in}

\begin{abstract}
Let $C(G)$ denote the set of lengths of cycles in a graph $G$. In
the first part of this paper, we study the minimum possible value of
$|C(G)|$ over all graphs $G$ of average degree $d$ and girth $g$.
Erd\H{o}s~\cite{Er1} conjectured that $|C(G)| =
\Omega\big(d^{\lfloor (g-1)/2\rfloor}\big)$ for all such graphs, and
we prove this conjecture. In particular, the longest cycle in a graph of
average degree $d$ and girth $g$ has length $\Omega\big(d^{\lfloor (g-1)/2\rfloor}\big)$.
The study of this problem was initiated by Ore in 1967 and
our result improves all previously known lower bounds on the length of the
longest cycle ~\cite{EM,EFRS,O,Z,ZB}.
Moreover, our bound cannot
be improved in general, since known constructions of $d$-regular
Moore Graphs of girth $g$ have roughly that many vertices.
We also show that $\Omega\big(d^{\lfloor (g-1)/2\rfloor}\big)$ is a lower bound
for the number of odd cycle lengths in a graph of chromatic number
$d$ and girth $g$. Further results are obtained for the number of
cycle lengths in $H$-free graphs of average degree $d$.
\medskip

In the second part of the paper, motivated by the conjecture of
Erd\H{o}s and Gy\'{a}rf\'{a}s~\cite{EG} (see also
Erd\H{o}s~\cite{E}) that every graph of minimum degree at least
three contains a cycle of length a power of two, we prove a general
theorem which gives an upper bound on the average degree of an
$n$-vertex graph with no cycle of even length in a prescribed
infinite sequence of integers. For many sequences, including the
powers of two, our theorem gives the upper bound $e^{O(\log\!^*n)}$
on the average degree of graph of order $n$ with no cycle of length
in the sequence, where $\log\!^*n$ is the number of times the binary
logarithm must be applied to $n$ to get a number which is at most
one.
\end{abstract}

\section{Introduction}

For a graph $G$, let $C(G)$ denote the set of integers whose
elements are lengths of cycles in $G$. The study of cycles in graphs
has long been fundamental, and many questions about properties of
graphs that guarantee some particular range of cycle length have
been considered. The central goal in this paper is to obtain lower
bound on $|C(G)|$ when $G$ is a graph of average degree $d$ and
girth $g$, or $G$ is an $H$-free graph, and to determine which
integers are guaranteed to appear in $C(G)$ under density conditions
on $G$.

\medskip

In the case of dense graphs, there are many results which determine
when $C(G)$ is an interval, or almost an interval. One of the first
results in this direction was obtained by Bondy~\cite{B}, who proved
that for every $n$-vertex graph $G$ of minimum degree larger than
$\frac{n}{2}$, $C(G) = \{3,4,,\dots,n\}$. Once the minimum degree of
a graph is allowed to pass below $\frac{n}{2}$, we can no longer
guarantee any odd cycles, as the graph may be bipartite. Also, one
cannot guarantee a hamiltonian cycle, (equivalently, $n \in C(G)$).
In this situation, the natural question is to ask when $C(G)$
contains all even integers up to $2\ell$, where $2\ell$ is the
length of a longest even cycle of $G$. Bollob\'{a}s and
Thomason~\cite{BT} showed that this is the case if $G$ has order $n$
and size at least $\lfloor\frac{n^2}{4}\rfloor - n + 59$. The best
result on this problem is by Gould, Haxell, and Scott~\cite{GHS},
who proved that if an $n$-vertex graph $G$ has minimum degree at
least $cn$, where $c > 0$ is a constant, then $C(G)$ contains all
even integers up to $2\ell - K$ for some constant $K$ depending only
on $c$. It is conjectured that for some constant $c > 0$, every
hamiltonian $n$ by $n$ bipartite graph of minimum degree at least
$c\sqrt{n}$ contains cycles of all even lengths in
$\{4,6,8,\dots,2n\}$. All of the above mentioned results are for
dense graphs -- graphs whose average degree is linear in the order
of the graph. In this paper, we are interested in studying $C(G)$
for sparse graphs.

\bigskip

Since $n$-vertex graphs of average degree $d$ may have girth at
least $\log_{d-1} n$, it is clear that for sparse graphs one cannot
hope to state that $C(G)$ contains any integer from a finite set.
Erd\H{o}s and Hajnal~\cite{EG} conjectured
\[ \sum_{\ell \in C(G)} \frac{1}{\ell} \; \; = \; \;  \Omega(\log d)\]
whenever $G$ has average degree $d$. (Here and throughout the paper the
notation $a_d=\Omega(b_d)$ means that there is an absolute constant $C$ such that
$a_d\geq Cb_d$ when $d \rightarrow \infty$.)
This conjecture was proved by
Gy\'{a}rf\'{a}s, Koml\'{o}s and Szemer\'{e}di~\cite{GKS}. Their
result shows that if a graph does not have too many short cycles,
then it must have many long cycles. However, it appears to be very
difficult to pass from such statement to statements about the size
or arithmetic structure of $C(G)$.

\subsection{Cycles in graphs of large girth}

The first problem we study is to determine the size of $C(G)$
for graphs of given average degree and girth (the {\sl girth} of $G$
is the length of the shortest cycle in $G$). The length of a longest
cycle in a graph of girth $g$ was first studied by Ore~\cite{O}.
For graphs of girth at most four and average degree $d$,
it is straightforward to prove that the longest cycle has length at
least $d + 1$ if the girth is three, and at least $2d$ if the girth is four,
and the proofs of these facts show $|C(G)| \geq d - 1$, with equality for $K_{d + 1}$ and $K_{d,d}$. Erd\H{o}s~\cite{Er1}
conjectured
that $|C(G)| = \Omega(d^{\lfloor (g - 1)/2 \rfloor})$ whenever
$G$ has average degree $d$ and girth $g$. This was proved for $g=5$ by
Erd\H{o}s, Faudree, Rousseau and Schelp~\cite{EFRS}. They also show that
$|C(G)| = \Omega(d^{5/2})$ for $g = 7$, $|C(G)| \geq \Omega(d^{3})$ for $g = 9$
and $|C(G)| =  \Omega(d^{g/8})$ in general. For comparison,
Erd\H{o}s' conjecture is $|C(G)| = \Omega(d^3)$ for $g = 7$ and
$|C(G)| = \Omega(d^4)$ for $g = 9$.
In Section 2, we will give a short proof of Erd\H{o}s' conjecture. In fact, we
obtain the following stronger theorem.

\begin{theorem}\label{main}
Let $G$ be a graph of average degree $d$ and girth
$g$. Then $C(G)$ contains $\Omega\big(d^{\lfloor (g - 1)/2 \rfloor}\big)$
consecutive even integers.
\end{theorem}

The study of the length of a longest cycle in graphs of girth $g$ and average or minimum degree
$d$ was initiated by Ore \cite {O} in 1967 and attracted attention of a number of
researchers~\cite{EM,EFRS,Z,ZB} since then. Our lower bound on $|C(G)|$
improves all of these results, and it best possible up to constant
factors: to see why Theorem \ref{main} cannot be improved, recall
that the {\sl Moore Bound} for a graph $G$ of minimum degree $d$ and
girth $g$ states
\begin{eqnarray*}
|V(G)| \geq \left\{\begin{array}{ll}
1 + d + d(d - 1) + \dots + d(d - 1)^{\lfloor \frac{g-1}{2}\rfloor -1}  & \mbox{ if }g\mbox{ is odd}\\
2\Big(1 + (d - 1) + (d - 1)^2 + \dots + (d - 1)^{\lfloor \frac{g-1}{2}
\rfloor}\Big) & \mbox{ if
}g\mbox{ is even}\\ \end{array}\right.
\end{eqnarray*}
Up to the constant factor, it is known that this bound is tight for
infinitely many values of $d$ whenever $g \leq 8$ or $g = 12$ and it
is also believed that for all other values of $g$ there are graph
with girth $g$  and order $O\big(d^{\lfloor (g-1)/2 \rfloor}\big)$.
So it is evident that $|C(G)| = O\big(d^{\lfloor (g-1)/2
\rfloor}\big)$ for such graphs. A related problem is to determine
the number of odd integers in $C(G)$ when $G$ has large chromatic
number and girth. For example, Gy\'{a}rf\'{a}s~\cite{G} proved that
a graph of chromatic number at least $2d + 1$ contains cycles of $d$
distinct odd lengths, and equality holds only for graphs all of whose
blocks are complete graphs. Using similar techniques as in the proof of Theorem \ref{main} we can
generalize the result of Gy\'{a}rf\'{a}s as follows: if $G$ is a
graph of chromatic number $d$ and girth $g$, then $C(G)$ contains
$\Omega\big(d^{\lfloor (g-1)/2 \rfloor}\big)$ consecutive integers.

\subsection{Cycles in $H$-free graphs}

A graph is {\sl $H$-free} if it contains no subgraph isomorphic to
$H$. We consider the following generalization of Theorem \ref{main}
in Section 3. Given a bipartite graph $H$, determine a lower bound
for $|C(G)|$ when $G$ is an $H$-free graph of average degree $d$.
Specifically, we consider $r$-half-bounded bipartite graphs. A
bipartite graph is {\sl $r$-half-bounded} if the degrees of all the
vertices in one color class are at most $r$. An example of such
graph is a complete bipartite graph $K_{r,s}$ with parts of size
$r\leq s$. Using recent estimates on Tur\'{a}n numbers for
$r$-half-bounded  graphs due to Alon, Krivelevich and
Sudakov~\cite{AKS}, we prove the following result.

\begin{theorem}\label{main2}
Let $H$ be a fixed bipartite graph containing a cycle and let $G$ be
an $H$-free graph of average degree $d$. Then there exists a
constant $t > 1$ depending on $H$ such that $C(G)$ contains
$\Omega\big(d^{t/(t - 1)}\big)$ consecutive even integers.
Furthermore, we can take $t = r$ if $H$ is $r$-half-bounded, and $t
= 1 + \frac{1}{k - 1}$ if $H$ is a $2k$-cycle.
\end{theorem}

Notice that when $H$ is a $2k$-cycle, this result generalizes our
Theorem \ref{main} from graphs of girth $2k + 1$ or $2k + 2$ to graphs with no $2k$-cycle.
The estimate for $r$-half-bounded graphs in Theorem \ref{main2} is
tight for every value of $r \geq 2$. Indeed, by the construction of
projective norm graphs in
\cite{ARS} (modifying that in \cite{KRS}) for every fixed
$s\geq (r-1)!+1$ there are graphs of order $O\big(d^{r/(r - 1)}\big)$ and
average degree $d$ which do not contain copy of
$K_{r,s}$.

\subsection{Arithmetic structure of $C(G)$}

In the second part of the paper, we discuss the arithmetic structure
of $C(G)$ for sparse graphs. The type of question we would like to answer is:
what is the smallest $d$ such that every graph of average degree at least $d$
has a cycle of length equal to a square, or a power of two, or twice a prime?
Our main theorem is motivated by the conjecture of Erd\H{o}s
and Gy\'{a}rf\'{a}s~\cite{EG}, stating that every graph of minimum degree
at least three contains a cycle of length a power of two, and by
the questions posed by Erd\H{o}s (see page 228 of~\cite{E}). Throughout
this section, for a sequence $(\sigma(i))_{i \geq 1}$ of integers, a
{\sl $\sigma$-cycle} is a cycle of length $\sigma(i)$ for some $i \geq 1$.
 We write $\pi < \sigma$
to denote that $\pi$ is a subsequence of $\sigma$.

\bigskip

Perhaps the most natural starting point is to determine when $C(G)$
contains an integer congruent to zero modulo a given integer $k$.
The first result in this direction was proved by
Bollob\'{a}s~\cite{Bol}, who showed that if $G$ has average degree
at least $\frac{2}{k}(k+1)^k$, then $G$ contains a cycle of length
zero modulo $k$. The main result in~\cite{V} (see also Fan~\cite{F})
shows that if $\sigma$ is any infinite increasing sequence of even
integers such that $|\sigma(j) - \sigma(j - 1)| \leq k$ for all $j
\geq 2$, then every graph of average degree at least $4k$ contains a
$\sigma$-cycle. In this section, we are interested in extending this
result to the case that $|\sigma(j) - \sigma(j - 1)|$ is not
bounded. The theorem below gives an upper bound on the average
degree of a graph containing no $\sigma$-cycles. In this theorem,
all logarithms are natural logarithms.

\begin{theorem}\label{main3}
For any infinite increasing sequence $\sigma$ of positive even
integers and for any $n$-vertex graph $G$, if $G$ contains no
$\sigma$-cycle, then $G$ has average degree at most
\[ \inf_{{\pi < \sigma}\atop{r \geq 1}}
\exp \left(6r + \sum_{i = 1}^{r}\frac{2\log \Delta(i)}{\pi(i - 1)} +
\frac{2\log n}{\pi(r)}\right),\] where $\pi(0) := 1$, $\Delta(1) :=
\pi(1)$, and $\Delta(i) = \max\{\sigma(j) - \sigma(j-1) : \sigma(j)
\leq \pi(i)\}$ for $i \geq 2$.
\end{theorem}

To illustrate this statement consider the case when $\sigma(i) = 2^i$
for $i \geq 1$. Then we can take $\pi$
to be the sequence of towers of twos, namely

\vspace{-0.2in}

\[ \pi(1) = 2 \hspace{0.2in} \pi(2) = 2^2 \hspace{0.2in} \pi(3) =
2^{2^{2}} \hspace{0.2in} \cdots 
%\hspace{0.2in} \pi(i) =2^{2^{2^{\iddots^{\;_{2}}}}}
\] so that $\pi(i) = 2^{\pi(i - 1)}$,
and take $r = \log\!\mbox{*}n$, where $\log\!\mbox{*}n = i$ whenever
$\pi(i - 1) < n \leq \pi(i)$.
  Then
Theorem \ref{main3} implies that every graph of order $n$ with no
cycle of length a power of two has average degree $\exp(O(\log^{*}
n))$. In fact, the same bound holds for many sequences, such as
twice primes, squares, and the tower sequence $\pi$ defined above.
We say that a sequence $\sigma$ is {\sl exponentially bounded} if
there is an absolute constant $C > 1$ such that $\sigma(i) \leq
C\sigma(i - 1)$ for all $i \geq 2$.

\begin{corollary}\label{exp}
Let $\sigma$ denote an infinite increasing exponentially bounded
sequence of positive even integers. Then any $n$-vertex graph with
no $\sigma$-cycles has average degree $\exp(O(\log\!\mbox{*}n))$.
\end{corollary}

This corollary will be proved in Section 4. Also in Section 4, we
will construct sequences $\sigma$ and $n$-vertex graphs with no
$\sigma$-cycles whose average degrees have the same order of
magnitude as the upper bound in Theorem \ref{main3}, up to absolute
constant factors in the exponent. However, these sequences are not
exponentially bounded. It would be interesting to see if there are
graphs with arbitrarily large average degree containing no
$\sigma$-cycles, for some exponentially bounded sequence $\sigma$ of
positive even integers.

\setcounter{section}{1}

\section{Cycles in graphs of large girth}

The (open) neighborhood of $X \subset V(G)$ in a graph $G$ is
defined by
\[ \partial X = \{y \in V(G) \backslash X \; | \;  \exists x \in X : \{x,y\} \in E(G)\}.\]
In words, this is the set of vertices not in $X$ and adjacent to at
least one vertex of $X$. The $d$-core of a graph $G$, when it
exists, is the subgraph obtained by repeatedly deleting vertices of
degree at most $d - 1$. It is a
well-known fact that if a graph has integer average degree $2d$,
then it has a $d$-core. It is convenient to assume throughout that
$d$ is an integer. Our first lemma states that graphs of large
average degree and girth expand on small sets.

\begin{lemma}\label{cyclexpand}
Let $G$ be a graph of girth $g$ and minimum degree at least $6(d +
1)$. Then, for every $X \subset V(G)$ of size at most
$\frac{1}{3}d^{\lfloor (g-1)/2 \rfloor}$, \[ |\partial X| >
2|X|.\]
\end{lemma}

\begin{proof}
Suppose $|\partial X| \leq 2|X|$ for some $X \subset V(G)$. Let $H$
be the subgraph of $G$ spanned by the set $Y=X \cup \partial X$.
Then $|Y|\leq 3|X|$ and \[ e(H)\geq \frac{1}{2}\sum_{x \in X} d(x)
\geq 3(d + 1)|X| \geq (d+1)|Y|.\] Thus $H$ contains a subgraph
$\Gamma$ with minimum degree $d+1$. Applying the Moore Bound to
$\Gamma$, we obtain:
$$ 3|X|\, \geq \,|Y| \, \geq \, |V(\Gamma)| \, \geq \, 1\, + \, (d + 1) \sum_{i < \lfloor
(g-1)/2 \rfloor} d^i\,
  >\,  d^{\lfloor (g-1)/2 \rfloor}$$
and therefore $|X| > \frac{1}{3}d^{\lfloor (g-1)/2 \rfloor}$, as required.
\end{proof}

\bigskip

Using Lemma \ref{cyclexpand}, we prove the conjecture of Erd\H{o}s
stating $|C(G)| = \Omega\big(d^{\lfloor (g-1)/2 \rfloor}\big)$ when $G$ has girth $g$ and
average degree $d$. A key ingredient of the proof is a lemma of
P\'{o}sa~\cite{P} (see also \cite{Lo}, Exercise 10.20) which says
that if $G$ is a graph and $|\partial X|
> 2|X|$ for every $X \subset V(G)$ of size at most $m$, then $G$
contains a path of length $3m$.

\begin{theorem}\label{numcyc}
For any graph $G$ of girth $g$ and average degree $48(d + 1)$,
$|C(G)| \geq \frac{1}{8}d^{\lfloor (g-1)/2 \rfloor}$.
\end{theorem}

\begin{proof}
Let $H$ be a maximum bipartite subgraph of $G$, containing at least half of the edges
of $G$. Then some connected component $F$ of $H$ has average degree at least
$24(d + 1)$. Let $T$ be a breadth first search tree in $F$, and let
$L_i$ denote the set of vertices of $T$ at distance $i$ from the
root of $T$. Since $F$ is bipartite, no edge of $F$ joins two vertices
of $L_i$. Denote by $e(L_i,L_{i+1})$ the number of edges of $F$ with
one endpoint in $L_i$ and one endpoint in $L_{i+1}$. Then
\begin{eqnarray*}
\sum_i e(L_i, L_{i+1}) &=& e(F)\, \geq\, 12(d+1)|V(F)|\, = \, 12(d+1)\sum_i |L_i|\\
&=&
6(d+1)\sum_i \big(|L_i|+|L_{i+1}|\big).
\end{eqnarray*}
Thus, there is an index $i$ such that the subgraph $F_i \subset F$
induced by $L_i \cup L_{i+1}$ has average degree at least $12(d +
1)$. Then $F_i$ contains a subgraph $\Gamma$ with minimum degree
$6(d + 1)$. By Lemma \ref{cyclexpand} we have that $|\partial X| >
2|X|$ for every $X \subset V(\Gamma)$ of size at most
$\frac{1}{3}d^{\lfloor (g-1)/2 \rfloor}$. Hence $\Gamma$ contains a
path $P$ of length $d^{\lfloor (g-1)/2 \rfloor}$ by P\'{o}sa's
Lemma. Let $T'$ be a minimal subtree of $T$ whose set of end
vertices is exactly $V(P) \cap L_i$. The minimality of $T'$ ensures
that it branches at the root. Let $A$ be the set of vertices in
$V(P) \cap L_i$ in one of these branches and let $B=(V(P) \cap
L_i)\setminus A$. Then both $A,B$ are nonempty and all paths from
$A$ to $B$ through the root of $T'$ have the same length, say $2h$.
We may assume that $|B|\geq |A|$ and $|B| \geq \frac{1}{4}|P|$. Let
$a$ be an arbitrary vertex in $A$. Then there are at least
$\frac{1}{2}|B| \geq \frac{1}{8}|P|$ vertices of $B$ on the same
side of path from $a$. Hence there  are subpaths of $P$ from $a$ to
a vertex of $B$ of at least $\frac{1}{8}|P|$ different lengths. For
any such path $Q$, there is a unique subpath $R$ of $T'$ through the
root joining the endpoints of $Q$, so that $Q \cup R$ is a cycle in
$G$. Since all $R$ have the same length $2h$, we obtain
$\frac{1}{8}d^{\lfloor (g-1)/2 \rfloor}$ cycles of different
lengths, and $|C(G)| \geq \frac{1}{8}d^{\lfloor (g-1)/2 \rfloor}$.
\end{proof}

\subsection{Proof of Theorem \ref{main}}

To obtain Theorem \ref{main}, we will slightly modify the proof of
Theorem \ref{numcyc}. A {\sl $\theta$-graph} is a graph consisting
of three internally disjoint paths between two vertices. We observe
the following lemma as a corollary of the proof of Theorem
\ref{numcyc}:

\begin{lemma}\label{theta}
Let $G$ be a graph of average degree $48(d + 1)$ and girth $g$,
where $d^{\lfloor (g - 1)/2 \rfloor} \geq 6$. Then $G$ contains a
$\theta$-graph containing a cycle of length at least $d^{\lfloor
(g-1)/2 \rfloor} + 2$.
\end{lemma}

\begin{proof}
Let the path $P$, tree $T'$ and set $L_i$ be defined as in the proof
of Theorem \ref{numcyc}. Since $d^{\lfloor (g - 1)/2 \rfloor} \geq
6$, we have $|V(P) \cap L_i| \geq 3$.
Let $Q \subset P$ be a path of length at least $|E(P)| - 2$ with endpoints in $L_i$.
Then also $|V(Q) \cap L_i| \geq 3$ and therefore $Q$ has an interior
vertex in $L_i$. If $R$ is a
path in $T'$ joining the endpoints of $Q$, then $Q \cup R$ is a
cycle of length at least $d^{\lfloor (g-1)/2 \rfloor} + 2$. Finally,
for some path $S \subset T'$ from the root of $T'$ to an interior
vertex of $Q$ in $L_i$, the subgraph $Q \cup R \cup S$ is the required
$\theta$-graph.
\end{proof}

\bigskip

It is convenient to define an {\sl $AB$-path} in a graph $G$ to be a
path with one endpoint in $A$ and one endpoint in $B$, where $A,B
\subset V(G)$. The following result of Bondy and
Simonovits~\cite{BS} (see also \cite{V}) will be used to prove Theorem \ref{main}.

\begin{lemma}\label{abpaths}
Let $\Gamma$ be a $\theta$-graph and let $(A,B)$ be
a nontrivial partition of $V(\Gamma)$.
Then $\Gamma$ contains $AB$-paths of all lengths less than $|V(\Gamma)|$
unless $\Gamma$ is bipartite with bipartition $(A,B)$.
\end{lemma}

{\bf Proof of Theorem \ref{main}.}
Let $G$ be a graph of average degree $192(d+1)$ and girth $g$
and let $H$ be a maximum bipartite subgraph of $G$. Then some connected component $F$
of $H$ has average degree at least $96(d+1)$.
Let $T$ be a breadth-first search tree in $F$, and let $L_i$ denote the set of vertices
of $T$ at distance $i$ from the root. Then, for some $i$, the subgraph $F_i$ of $F$
induced by $L_i \cup L_{i+1}$ has average degree at least $48(d +
1)$. By Lemma \ref{theta}, $F_i$ contains a $\theta$-graph $\Gamma$
containing a cycle of length at least $d^{\lfloor (g-1)/2 \rfloor} + 2$. Let $T'$ be the
minimal subtree of $T$ whose set of end vertices is $V(\Gamma) \cap L_i$.
Then there is a partition $(A,B^*)$ of $V(\Gamma) \cap L_i$ such that all
$AB^*$-paths in $T'$ go through the root and have the same length, say $2h$.
Let $B = V(\Gamma) \backslash A$. By Lemma \ref{abpaths}, there exist $AB$-paths in $\Gamma$
of all even lengths in $\{1,2,\dots,d^{\lfloor (g-1)/2 \rfloor} + 2\}$. Since they have an
even length, each such path
is actually an $AB^*$-path, and the union of this path with the unique subpath of
$T'$ of length $2h$ joining its endpoints is a cycle. Therefore
$C(G)$ contains $d^{\lfloor (g-1)/2 \rfloor}$ consecutive even integers, as required.
\hfill{ }\vrule height10pt width5pt depth1pt

\subsection{Chromatic number and cycle lengths}

Using the above methods, we prove that in a graph $G$ of large
chromatic number and girth, $C(G)$ contains long interval of
consecutive integers. We only sketch the details, as they resemble
the proof of Theorem \ref{main}. First we require two simple
lemmas.

\begin{lemma}\label{oddeven}
Let $H$ be a minimal $d$-chromatic graph, where $d \geq 3$. Then for any
distinct vertices $u,v \in V(H)$, there is a $uv$-path of odd length in $H$ and
a $uv$-path of even length in $H$.
\end{lemma}

\begin{proof}
Since $H$ is minimal $d$-chromatic, and $d \geq 3$, $H$ has no cut-vertex.
Fix $u,v \in V(H)$, and an odd cycle $C \subset H$. By Menger's Theorem,
there exist two vertex-disjoint paths $P,Q$ starting at
$u,v$ and ending at vertices $w,x \in V(C)$, respectively.
Since $C$ is an odd cycle, $C = R \cup S$ where $R$ and $S$ are
internally disjoint $wx$-paths whose lengths have different parity.
It follows that $P \cup Q \cup R$ and $P \cup Q \cup S$ are $uv$-paths whose lengths have
different parity.
\end{proof}

\bigskip

A $\theta$-graph is odd if it is non-bipartite. 

\begin{lemma}\label{oddtheta}
Let $H$ be a minimal $d$-chromatic graph, where $d \geq 3$. Then for any
even cycle $C \subset H$, there is an odd $\theta$-graph in $H$ containing $C$. 
\end{lemma}

\begin{proof}
For $u,v \in V(C)$, let $d(u,v)$ be the distance from $u$ to $v$ on $C$.
By Lemma \ref{oddeven}, for any $u,v \in V(C)$ we can find a path $P$
such that $|E(P)| \not = d(u, v)~ (\hspace{-0.25cm}\mod 2)$. Let $P = (u_0,u_1,\dots,u_r)$ be the shortest path with $u_0, u_r \in V(C)$ and 
$|E(P)| \not = d(u_0,u_r) ~ (\hspace{-0.25cm}\mod 2)$. Let $Q \subset P$ be the path $(u_0,u_1,\dots,u_s)$ with $u_s \in V(C)$ and 
$u_i \not \in V(C)$ for $i < s$. 
If $|E(Q)|= d(u_0,u_s) ~ (\hspace{-0.25cm}\mod 2)$, then $R = P - \{u_i : i < s\}$ is a $u_s u_r$-path
with $|E(R)| \not = d(u_{s-1},u_r) ~ (\hspace{-0.25cm}\mod 2)$, contradicting the choice of $P$. 
So $|E(Q)| \not = d(u_0,u_s) ~ (\hspace{-0.25cm}\mod 2)$, and $C \cup Q$ is an odd $\theta$-graph.
\end{proof}

\bigskip

We now prove our main result concerning the number of odd cycle lengths for graphs of large chromatic number and girth:

\begin{theorem}\label{chi}
Let $G$ be a graph of chromatic number $d$ and girth $g$. Then
$C(G)$ contains $\Omega\big(d^{\lfloor (g-1)/2 \rfloor}\big)$
consecutive integers.
\end{theorem}

\begin{proof}
Take a breadth first search tree $T$ in a component of $G$ of
chromatic number $d$, and let $L_i$ denote the set of vertices at
distance $i$ from the root of $T$. Then, for some $i$, the subgraph
$F$ spanned by $L_i$ has chromatic number at least
$\frac{1}{2}d$. Let $H$ be a minimal $\frac{1}{2}d$-chromatic subgraph
of $G[L_i]$. By Lemma \ref{theta}, assuming $d$ is large, $H$ contains a $\theta$-graph containing
a cycle $C$ of length $\Omega\big(d^{\lfloor (g-1)/2 \rfloor}\big)$. 
By Lemma \ref{oddtheta}, we can ensure that this $\theta$-graph is an odd $\theta$-graph
containing $C$, which we denote by $\Gamma$. 
If $T'$ is a minimal subtree of $T$ whose set of end vertices is
$V(\Gamma)$, then $T'$ branches at the root of $T'$, and this gives
a partition $(A,B)$ of $V(\Gamma)$ as in the proof of Theorem
\ref{main}. By Lemma \ref{abpaths}, there are $AB$-paths of all
lengths less than $|V(\Gamma)|$, and these paths together with
subpaths of $T'$ give the required cycle lengths.
\end{proof}

\section{Cycles in $H$-free graphs}

To prove Theorems \ref{main2} and \ref{main3}, we will use the
following lemma, which summarizes the ideas of Section 2.
Recall that a property of graphs is {\sl monotone} if it is closed
under taking subgraphs. The {\sl radius} of graph $G$ is the smallest integer
$r$ for which there is a vertex $v$ in $G$ such that the distance from any other vertex
of $G$ to $v$ is at most $r$.

\begin{lemma}\label{largetheta}
Let ${\cal P}$ be a monotone property of graphs, and suppose that
for every graph $G \in {\cal P}$ with minimum degree $d$, and every set $X \subset
V(G)$ of size at most $f(d)$,
\[ |\partial X| > 2|X|.\]
Then every $G \in {\cal P}$ of average degree at least $16d$
contains cycles of $3f(d)$ consecutive even lengths, the shortest
having length at most twice the largest radius of any component of
$G$.
\end{lemma}

\begin{proof}
Let $G'$ be a maximum bipartite subgraph of $G$, and let $T$ be a breadth-first
search tree in a connected component $F$ of $G'$ of average degree at least
$8d$. If $L_i$ is the set of vertices at distance $i$ from the root
of $T$ in $F$, then for some $i$, the subgraph $F^*$ of $F$ induced
by $L_i \cup L_{i+1}$ has average degree at least $4d$. Let $T^*$ be
a breadth first search tree in a connected component of $F^*$ of average
degree at least $4d$. If $L_j^*$ is the set of vertices at distance
$j$ from the root of $T^*$ in $F^*$, then for some $j$, the subgraph
$F^*_j$ of $F^*$ induced by $L_j^* \cup L_{j+1}^*$ has average
degree at least $2d$. Now let $\Gamma$ be a subgraph of $F^*_j$ with minimum degree
at least $d$. Since ${\cal P}$ is monotone property and  $G\in {\cal P}$
we have that also $\Gamma \in {\cal P}$. Therefore, $|\partial X| > 2|X|$ for
every subset $X\subset V(\Gamma)$ of size at most $f(d)$.
By P\'{o}sa's Lemma, there is a path $P \subset \Gamma$ of length
$3f(d)$. If $T'$ is a minimal subtree of $T^*$ whose set of
end vertices is $V(P) \cap L_j^*$, then as in Lemma \ref{theta}, $P
\cup T'$ contains a $\theta$-graph, $J$, containing a cycle of
length $3f(d) + 2$. Let $T''$ be the minimal subtree of $T$ whose set
of end vertices is $V(J)$.
Applying Lemma \ref{abpaths} as in the proof
of Theorem \ref{main}, we see that $J \cup T''$ contains cycles of
$3f(d)$ consecutive even lengths in $G$. Since the shortest cycle
has length at most $2i + 2$, the proof is complete.
\end{proof}

\bigskip

In what follows, we denote by $\mbox{ex}(n,H)$ the {\sl Tur\'{a}n
number} of graph $H$, which is the maximum number of edges in an
$H$-free graph on $n$ vertices. To obtain expansion in $H$-free
graphs, where $H$ is bipartite, we show that it is enough to find
upper bounds for $\mbox{ex}(n,H)$.

\begin{lemma}\label{expand}
Let $a > 0, \frac{1}{2} < b < 1$ be reals such that for any
positive integer $n$, $\mbox{ex}(n,H) \leq a n^{2b}$.
Then, for any $H$-free graph $G$ of minimum degree at least
$18ad$, and any subset $X$ of vertices of $G$ of size at most
$d^{1/(2b - 1)}$, $|\partial X| > 2|X|$.
\end{lemma}

\begin{proof}
Suppose that $X$ is a subset of $G$ of size $m$ such that $|\partial X| \leq 2|X|$.
Let $G_Y$ be the subgraph of $G$ induced by $Y=X \cup \partial X$.
Then $|Y| \leq 3m$ and $e(G_Y) \geq 9ad|X|=9adm$. On the other hand, since $G_Y$ is $H$-free we have that
$$9adm\, \leq\, e(G_Y)\, \leq \, \mbox{ex}(|Y|,H)\, \leq\, a|Y|^{2b}\, \leq\,
a(3m)^{2b}\,<\,9am^{2b}.$$
Therefore $m^{2b - 1} > d$, which proves the lemma.
\end{proof}

\bigskip
\medskip

{\bf Proof of Theorem \ref{main2}.}
By the well known result of K\"ovari, S\'os and Tur\'an~\cite{KST},
for every bipartite graph $H$ there are two constants
$t > 1$ and $c$ depending only $H$ such that
$$ \mbox{ex}(n,H) \leq cn^{2 - 1/t}.$$
By Lemma \ref{expand}, with $a = c$ and $b = 1 - 1/2t$, every
$H$-free graph $F$ of minimum degree at least $18cd$ has the
property that for every $X \subset V(F)$ of size at most $f(d) =
d^{t/(t - 1)}$, $|\partial X| > 2|X|$. By Lemma \ref{largetheta},
with ${\cal P}$ equal to the set of all $H$-free graphs, we deduce that
every $G \in {\cal P}$ of average degree $288cd$ contains
$3f(d)$ cycles of consecutive even lengths, proving the theorem.
For the particular case when $H$ is $r$-half-bounded,
Alon, Krivelevich and Sudakov~\cite{AKS} showed that
$$ \mbox{ex}(n,H) = O(n^{2 - \frac{1}{r}}),$$
so the proof above applies with $b = 1 - 1/(2r)$ and gives
$\Omega\big(d^{r/(r-1)}\big)$ cycles of consecutive even lengths. Finally,
if $H = C_{2k}$, then Corollary 9 in~\cite{V} shows
$$\mbox{ex}(n,H) = 8k n^{1 + \frac{1}{k}},$$
so we can apply the above proof with $b = \frac{1}{2} + \frac{1}{2k}$
to conclude that for every $C_{2k}$-free graph $G$, $C(G)$ contains
$\Omega(d^k)$ consecutive even integers.
\hfill{ }\vrule height10pt width5pt depth1pt

\bigskip

\section{Arithmetic structure of $C(G)$.}

Fix $\pi < \sigma$, and let ${\cal P}_i, i \geq 1$ denote the monotone
property of graphs containing no cycle of length $\sigma(j)$ for all
$\sigma(j) \leq \pi(i)$, and recall
$\Delta(i) = \max\{\sigma(j) - \sigma(j - 1) : \sigma(j) \leq \pi(i)\}$.
To prove Theorem \ref{main3}, we first prove the
following claim.

\begin{claim}
\label{arithmetic} Let $(a_i)_{i \geq 1}$ be positive real numbers
such that $a_1 = 4\pi(1)$ and, for all $i \geq 2$,
$$ \pi(i - 1) \log \frac{a_i}{288a_{i-1}} \,\geq\, 2\log \Delta(i).$$
Then, for every $n$-vertex graph $G \in {\cal P}_i$,
$$e(G) \leq a_i n^{1 + \frac{2}{\pi(i)}}.$$
\end{claim}

\begin{proof}
We proceed by induction
on $i$. For $i = 1$, Corollary 9 in~\cite{V} gives
$$e(G) \leq 4\pi(1) n^{1 + \frac{2}{\pi(1)}},$$
and this proves the claim for $i = 1$. Suppose we have proved the
claim for $j < i$, and let $G \in {\cal P}_i$ be an $n$-vertex graph
with $e(G) > a_i n^{1+ 2/\pi(i)}$, where $a_i$ satisfies the bounds
in the claim. By the induction hypothesis we have that every
$m$-vertex graph in ${\cal P}_{i-1}$ has at most $a_{i-1} m^{1+
2/\pi(i-1)}$ edges. Therefore by Lemma \ref{expand}, we have that
for any graph $F \in {\cal P}_{i-1}$ with minimum degree $d$ every
subset $X\subset V(F)$ of size at most
\[ f(d) = \Bigl(\frac{d}{18a_{i-1}}\Bigr)^{\frac{1}{2}\pi(i - 1)}\]
has $|\partial X| > 2|X|$.

\medskip

Since $G$ has $n$ vertices and $e(G) \geq a_i n^{1+2/\pi(i)}$,
by Lemma 6 in~\cite{V}, there is a subgraph $\Gamma$ of $G$ of average degree at least $a_i$
and radius at most $\frac{1}{2}\pi(i)$. Note that $\Gamma$ has property ${\cal P}_i$ and
thus has also property ${\cal P}_{i-1}$. By Lemma \ref{largetheta}, $\Gamma$
contains cycles of at least $3f(\frac{a_i}{16})$ consecutive even
lengths, the shortest of which has length at most $\pi(i)$. Since $\Gamma
\in {\cal P}_i$, there must be less than $\Delta(i)$ of these
consecutive even lengths, otherwise $\Gamma$ contains a cycle of length
$\sigma(j)$ for some $j \leq \pi(i)$. Therefore
\[ \left(\frac{a_i}{288a_{i-1}}\right)^{\frac{\pi(i - 1)}{2}}
\; = \;f\Bigr(\frac{a_i}{16}\Bigl)\;<   \; 3f\Bigr(\frac{a_i}{16}\Bigl) \; < \; \Delta(i)\] which
contradicts the bounds on $a_i$ in the claim.
\end{proof}

\bigskip

{\bf Proof of Theorem \ref{main3}.}\,
Recall that $\pi(0)=1$, $\Delta(1)=\pi(1)$ and let
\[ a_r = 4\pi(1) (288)^{r-1}\prod_{i = 2}^{r} \exp\left(\frac{2\log \Delta(i)}{\pi(i - 1)}
\right) < \frac{1}{2}\exp\left(6r +
\sum_{i=1}^{r}\frac{2\log\Delta(i)}{\pi(i - 1)}\right).\] Since
$a_r$ satisfies the condition of the Claim \ref{arithmetic}, we have
that the estimate on the number of edges of $G$ from this claim is
valid for any $r$. Therefore the average degree of $G$ is at most
\[ \inf_{r \geq 1} 2a_r n^{2/\pi(r)} \, \leq  \,   \inf_{r \geq 1}
\exp\left(6r + \sum_{i=1}^{r}\frac{2\log \Delta(i)}{\pi(i - 1)} +
\frac{2\log n}{\pi(r)}\right).\] This bound is valid for any $\pi <
\sigma$, so this completes the proof of Theorem \ref{main3}. \hfill{
}\vrule height10pt width5pt depth1pt

\bigskip

{\bf Proof of Corollary \ref{exp}.}\, Since $\sigma$ is
exponentially bounded, $\sigma(i) \leq C\sigma(i - 1)$ for all $i
\geq 2$. Let $r = \log\!\mbox{*}n$ and let $\pi < \sigma$ be chosen
so that $2^{\pi(i - 1)} \leq \pi(i) \leq (2C)^{\pi(i - 1)}$ for all
$i \geq 2$. Note that $\pi(r) \geq n$. Then, since $\Delta(i) \leq
\pi(i)$, the upper bound in Theorem \ref{main3} is

\[ \exp\left(6r + \sum_{i=1}^{r} \frac{2\log \pi(i)}{\pi(i
- 1)} + \frac{2\log n}{n}\right) \leq \exp(6r  + 2r\log(2C) + 2) =
e^{O(\log\!\mbox{*}n)}.\]

This proves Corollary \ref{exp}. \hfill{ }\vrule height10pt width5pt
depth1pt

\bigskip

In conclusion, we show that the bound in Theorem \ref{main3} cannot
be improved in general.

\bigskip

{\bf Construction.}  We construct a sequence of graphs $G_{1},G_{2},G_{3},\dots $
such that $|V(G_k)| = n_k - 2$, $G_{k}$ is $(\alpha_k + 1)$-regular,
and $G_{k}$ has girth larger than $n_{k-1}$, where $n_k$ is even
for all $k$ and $n_0 := 2$.
By known probabilistic and explicit constructions of small graphs of large girth
(see, e.g., \cite{LPS} and \cite{M}), we may
take $\log n_k = n_{k-1}\log \alpha_{k}$. Now let $\sigma$ be defined
by $\sigma(i) = n_i$. Since $G_{i}$ has girth larger than $n_{i-1}$
and $G_{i}$ has order less than $n_{i}$, none of the graphs $G_i$
have a $\sigma$-cycle. We choose $\alpha_i$ so that
\[ \alpha_i \geq 2^{2^{2^i}}.\]
If we take $\pi = \sigma$ in Theorem \ref{main3}, and $r = i$, then
we have $\Delta(i) \sim n_i$ as $i \rightarrow \infty$, and also
\[ \frac{\log \alpha_{i}}{\sum_{j < i} \log \alpha_{j}} \rightarrow \infty.\]
Also, as $j \rightarrow \infty$,
\[ \frac{2\log\Delta(j)}{\pi(j-1)} \sim \frac{2\log
n_j}{n_{j-1}} = 2\log\alpha_{j},\] and the upper bound on the
average degree of $G_i$ from Theorem \ref{main3} is:
\begin{eqnarray*}
\exp\left(6i + \sum_{j=1}^{i} \frac{2\log \Delta(j)}{\pi(j-1)} +
\frac{2\log n_i}{\pi(i)}\right) &=& \alpha_i^{2 + o(1)}.
\end{eqnarray*}
Since $G_i$ has average degree $\alpha_i + 1$, the bound given by
Theorem \ref{main3} is tight up to the constant factor in the
exponent.

\section{Concluding Remarks}

\bigskip

$\bullet$ {\bf Even cycles.} It would be interesting to determine if
there is an infinite increasing exponentially bounded sequence
$\sigma$ for which there are graphs of arbitrarily large average
degree containing no $\sigma$-cycles. Erd\H{o}s~\cite{E} states that
this is probably true when $\sigma$ is the sequence of powers of
two, although no example of a graph of minimum degree three with no
cycle of length a power of two is known. The construction in Section
4 shows that if we take a sequence $\sigma$ of positive integers
defined by $\log \sigma(i) = 2^{2^{i}}\sigma(i - 1)$, for $i \geq
1$, then there are graphs of arbitrarily large average degree with
no $\sigma$-cycles.

\bigskip

$\bullet$ {\bf Odd cycles.} Erd\H{o}s~\cite{E} posed analogous questions for
odd cycles in graphs of large chromatic number, for example, does every graph of
infinite chromatic number contain a cycle of length equal to an odd
integer square? By repeating a similar construction to that given in
Section 4, it is possible to show that there are (very fast-growing)
infinite increasing sequences of odd integers $\sigma$
and graphs of infinite chromatic number containing no $\sigma$-cycle.
On the other hand, we ask the following concrete question:
does every graph of chromatic number at least four contain a cycle of length one more than a power of two? This seems to be a natural generalization
of the Erd\H{o}s-Gy\'{a}rf\'{a}s~\cite{EG} conjecture.

\bigskip

$\bullet$ {\bf Chromatic number and cycle lengths.} It seems that
the result of Theorem \ref{chi} can be further improved. In
particular in \cite{Er1}, Erd\H{o}s asked whether for every
$\epsilon > 0$ and sufficiently large $d$, every triangle free graph
of chromatic number $d$ contains at least $\Omega(d^{2 -\epsilon})$
cycles of different lengths. More generally one can ask if every
graph of girth at least $2t$ and chromatic number $d$ contains at
least $\Omega(d^{t -\epsilon})$ cycles of different length. We
believe that the techniques which were developed in this paper may
be useful to attack these problems.

\vspace{0.2cm}
\noindent {\bf Acknowledgment.}\, The authors would like to thank 
the referee for careful reading of this manuscript.

\end{document}